\documentclass[12pt]{elsarticle}
\usepackage[utf8]{inputenc}
\usepackage{amsmath}
\usepackage{amssymb}
\usepackage{graphicx}
\usepackage{tikz}
\usetikzlibrary{arrows,petri,topaths}%
\usepackage{tkz-berge}
\usepackage{soul}

\usepackage{subcaption}
\usepackage{mwe}

\newtheorem{example}{Example}
\definecolor{cblue}{rgb}{0,0,1.0}
\title{Controlling inventory on electric roads}
\author{Alejandro Gutierrez-Alcoba, Roberto Rossi, Belen Martin-Barragan, Tim Embley}

\begin{document}
	\begin{abstract}
	Electric road systems (ERS) are roads that allow compatible vehicles to be powered by grid electricity while in transit, reducing the need for stopping to recharge electric batteries. We investigate how this technology can affect routing and delivery decisions for hybrid heavy good vehicles (HGVs) travelling on a ERS network to support the demand of a single product faced by a set of retailers in the network. We introduce the Electric Roads Routing Problem, which accounts for the costs of electricity and fuel on a ERS network, consumption that are affected by the battery level of the vehicle in each step of the journey, the routing decisions and the variable weight of the vehicle, which depends on vehicle load and delivery decisions. In particular, we study a stochastic demand version of the problem, formulating a mathematical programming heuristic and proving its effectiveness. We use our model on a realistic instance of the problem, showcasing the different strategies that a vehicle may follow depending on fuel costs in relation to the costs of electricity.
	\end{abstract}
	
	\maketitle
	
	\section{Introduction}
	\label{sec:intro}
	In recent years, growing concerns on climate change has seen the introduction of new CO2 emissions targets in many countries. For instance, in the UK, the Government recently established a net-zero emissions target for 2050, in an amendment of the Climate Change Act (2008) \cite{UKtopia}. For the transportation sector, one of the largest contributors of greenhouse gas (GHG) emissions, this means a complete transformation is needed in coming years. The adoption of electric vehicles (EVs), coupled with increasing levels of renewable energy production, is seen as a solution without compromising future transportation demands. 
	
	However, limitations associated with EVs, such as battery range, are tighter for heavy good vehicles (HGV). HGV are, at the same time, the vehicles causing the greatest environmental impact: in the EU, freight transport by HGVs accounts for about a quarter of all CO2 emissions caused by road transport \cite{Schroten2012}. Still, some engineering solutions could see HGVs harnessing electric power without relying on large electric batteries that compromise the payload and vehicle's weight and put stress on natural resources needed to manufacture them. These solutions, named as Electric Roads, eRoads, eHighways or Electric Road Systems (ERS) allow compatible vehicles to use power from the electric grid to travel along, reducing or even eliminating the need to stop for recharging. These systems can be classified into three categories, including: overhead catenary cables; conductive power from road; inductive power from road.
	
	In recent years, overhead catenary cables systems have been trialled for HGV use in countries such as Germany \cite{ScaniaR450}. In the UK, the Centre for Sustainable Road Freight has proposed a plan \cite{Cebon2020} for decarbonising road freight transport in the country. It calls for covering about 65\% of the roads travelled by HGVs (7500 km of road) with catenary cables by the late 2030s, with a construction investment estimated at £19.3 billion. HGV vehicles would be equipped with a pantograph system that allow them to connect and disconnect from installed catenary cables during transit, without the need for stopping. While connected, the cables power the vehicle and can charge an electric battery that is used to power the vehicle on conventional roads. To ensure a HGV can continue travelling when the ERS network and the battery is not enough, especially during the transition period, an extra propulsion system (such as a diesel engine) can also be included. An example of such vehicle is the hybrid Scania R 450, used in the first trials in Germany \cite{ScaniaR450}.
	
	The decarbonisation of road freight transport by deploying ERS on road networks and using hybrid vehicles poses many interesting research questions in transportation, ranging from routing to network design. In this paper, we explore some of the challenges that emerge from this concept from the perspective of inventory control. Our interest is to study a lot sizing problem where the trade-off costs of fulfilling the demand of a single product are given in terms of the environmental impact of delivery, in the context of HGVs navigating ERS networks. We present our work as a new problem termed the Electric Roads Routing Problem (ERRP). Three types on inter-related inventories are considered in this problem: First, the inventory at a set of retailers aiming to satisfy the demand of certain product. Second, the inventory of such a product on the vehicle delivering it to the retailers. Finally, the battery energy inventory needed to deliver the product to the retailers through a network with on-road charging opportunities. The inventories are inter-related. Obviously, the decisions of delivering product to certain retailer affect the vehicle product load, and consequently the availability of products for other retailers, or the need of the vehicle to visit the depot to load new product. In addition, we also consider that the energy consumption depends on the weight of the vehicle, meaning that the energy inventory is also affected by the product load of the vehicle.
		
	Our article also contributes to the literature on energy and pollution based models in vehicle routing, since we include important features that are relevant for hybrid vehicles that can take advantage of ERS networks. The problem considers the energy requirements of the vehicle and features recharging decisions, which occurs when a vehicle transits an ERS segment of road, and delivery decisions to retailers, considering the impact of weight on energy requirements. In particular, we study a stochastic demand setting of the problem and develop mathematical programming heuristics to achieve nearly optimal solutions.

	
	We contribute to the literature as follows:

	\begin{itemize}
	    \item We introduce the Electric Roads Routing Problem (ERRP), effectively accounting for battery and fuel requirements of a vehicle while transiting an ERS network to fulfil the demand of retailers.
	    \item We develop a stochastic dynamic programming (SDP) formulation to tackle the ERRP in a general setting where demand is stochastic.
	    \item We introduce a mixed integer linear programming (MILP) heuristic to tackle the stochastic variant of the ERRP, and demostrate the effectiveness of the proposed method on an extensive test bed, contrasting the solutions obtained by the heuristic against the optimal solutions obtained via SDP.
	    \item We test our heuristic on a realistic case study set in the area of Kent (UK), and perform a sensitivity analysis on fuel cost.
	\end{itemize}
	
	The rest of this document is structured as follows. In Section \ref{sec:litsurvay} we survey literature on inventory and routing problems and position our contribution. Section \ref{sec:ERRP} describes the problem setting of the ERRP. In Section \ref{sec:SDP-SERRP} we introduce a stochastic demand variant and present an SDP formulation. Section \ref{sec:MILPformulation} presents the stochastic MILP heuristic. A computational study and results proving the effectiveness of the heuristics are provided in Section \ref{sec:comp_study}. In section \ref{sec:kent} we present a case study of the problem and investigate the effect of the costs of fuel on the solutions found by our model. Finally, our conclusions are drawn in Section \ref{sec:conclusions}.
	
	\section{Literature survey}
	\label{sec:litsurvay}
	
	Our work can be positioned within two main strands of literature: inventory control and routing. Our starting point is of a set of retailers with their own inventory and periodic demand of a single product, with a replenishment plan that depends on the route followed by a single vehicle to deliver the product. In this section we review relevant literature of both strands, pointing out the differences with our approach.
	
	Regarding inventory, this paper deals with a non-stationary stochastic lot sizing problem: it is a single item inventory control problem, defined over a finite time horizon, with stochastic non-stationary demand. The literature in inventory control with stationary demand is vast, but academic research in lot sizing with non-stationary stochastic demand is not. Dynamic lot sizing was first introduced in \cite{WagnerWhithin}, where an exact solution algorithm is given. In \cite{Scarf1960}, the optimality of (s,S)-type policies for the lot sizing problem under stochastic demand is demonstrated. \cite{Bookbinder1988} presents a framework consisting of three control strategies for stochastic demand lot sizing problems: static uncertainty, static-dynamic uncertainty and dynamic uncertainty. In the first one, timing and quantity of orders are fixed at the beginning of the planning horizon. In the static-dynamic uncertainty strategy, order timings are fixed at the beginning of the time horizon, but the quantities are decided at the time of the replenishment, when the current inventory is known. In the third and most flexible strategy, the decision maker observes the inventory at the beginning of each period, deciding in first place if an order should be placed, and the replenishment quantity in case affirmative. This work motivated further studies aimed at determining optimal replenishment plans under these strategies. For instance \cite{SOX1997} and \cite{Vargas2009} focus on the static strategy. For the static-dynamic strategy, \cite{TARIM2006} presents a MILP formulation, modeling service quality under penalty costs. More recently, \cite{ROSSI2015} provides a model based on piecewise linearisation of the first order loss function and of its complementary. Our work  departs from the lot sizing problem literature in two main aspects. First, replenishment timings are constrained by the modeled vehicle supporting retailers: a replenishment is only possible during periods at which the vehicle is situated at a retailer position. Secondly, lot sizing problems model set up and product quantity costs typically using the linear function $k + Qc$, where $k$ is the set up order cost, $Q$ is the quantity and $c$ the unit procurement cost. In this work, delivering product to a retailer incurs costs in terms of the energy consumption of the vehicle, which takes into account the weight of the vehicle in each period. In this way, our problem considers both unit procurement costs and set up costs using a realistic physics model (see \ref{sec:energymodel}), which also depends on road conditions and the electric battery level at any period.
	

	There are studies that aim to incorporate a routing element to the problems as well. The inventory routing problem (IRP) and its variants \cite{IRPsReview} is perhaps the most well studied problem that combines both inventory and routing elements. In its basic form, the IRP consist of a set of retailers holding inventory of a single product that is depleted on each period by the demand they face. There is a single warehouse, from which each vehicle starts and finishes a route on any single period, visiting retailers along the route. The ERRP shares many features with the IRP, but the key difference is that in each time period the vehicle moves from one node of the road network to an adjacent one, rather than performing a complete route. In this sense, we present the ERRP as a \emph{small bucket} model \cite{Belvaux2001993}, in which time intervals are small, treating in detail what happens during them. 

	Our route modeling strategy not only differentiate our approach from this literature but from the routing literature in general. Since its introduction in \cite{Dantzig-truck}, the Vehicle Routing Problem (VRP) has been a central framework for routing research, see e.g. \cite{reviewVRP1}. In its basic version, the problem consists of devising a least-cost set of routes on a given network, visiting a set of retailers and with all routes starting and finishing at a single location. Many variants have been studied including other side constraints, such as the capacitated VRP (CVRP) \cite{CVRPalgorithms} or the VRP with Time Windows (VRPTW) \cite{VRPTWalgorithms}. 
	
	Common objective functions for the VRP minimise magnitudes that are proportional to distance. The Pollution Routing Problem (PRP) \cite{PRP-2011}, based on the VRP, sets a different approach for the objective by incorporating a comprehensive model, adopted from \cite{Barth2004}, to estimate energy consumption and GHG emissions of the vehicles. The model considers vehicle parameters related to forces against motion, such as aerodynamics and rolling resistance, and vehicle parameters and decision variables such as cruising speed and the weight of the vehicle to obtain an estimation of the required power in each step of the route. Since then, this approach has been followed in other green logistic studies that are not directly extensions, or variants, of the PRP (e.g., \cite{Goeke2015}, \cite{Lu2020}).
	
	Another important aspect in green logistics is on recharging decisions on route. The Green-VRP \cite{Erdogan2012} is introduced as a minimising environmental impact problem for scheduling routes of vehicles with a limited refueling infrastructure: in addition to the required retailer visits, the vehicles may visit refueling nodes when needed. Motivated by this problem, the authors in \cite{Schneider2014} introduce the Electric-VRPTW (E-VRTPTW), considering as well the time spent for complete battery recharges at stations for EVs and including time windows. Later in \cite{Goeke2015}, the authors extend the E-VRTPTW including a realistic model for energy consumption (as in \cite{PRP-2011}), to determine more accurately recharging times. Other relevant studies on recharging decisions for EVs are \cite{Schiffer2017}, in which partial recharges are allowed or \cite{Keskin2019}, which introduces the idea of time-dependent queuing times at the charging stations. 

    Table \ref{table:literature} classifies the surveyed literature in terms of the main features that differentiate our work.
    As in the E-VRPTW and related literature, charging decisions are of key importance in the ERRP, but in our case charging opportunities take place when a vehicle transits an ERS road segment, rather than stopping at certain nodes (stations). This is also the case in some previous studies (\cite{Artmeier-LNAI2010,Eisner2011,Storandt2012}) consisting of finding the best route from an origin to a destination for an EV, albeit the possibility of recharging is only on descending road segments, where regenerative braking can be applied. This sets an important difference: energy recovery is subject to the law of conservation of energy and, unlike in the ERRP, navigating in a close loop will always render less available battery at the end of the loop. Another important difference with this literature is that the ERRP is concerned with serving a set of retailers along the route, rather than on finding a shortest path between a pair of origin and destination nodes. To the best of our knowledge, there are no research studies on green logistic problems featuring recharging decisions on ERS networks.
    
	Classical transportation problems such as the VRP and their variants are mathematically defined on a \emph{complete graph}, with the abstraction of each arc between two retailers representing the shortest path of a sparsely connected graph representing the \emph{original road graph} \cite{Toth-Vigo}. This is the common case in the literature: the cost weights of the arcs represent distance or any other figure (travel time, energy consumption, fuel, fuel costs, emissions) that, in the context of the problem, is highly correlated with distance. However, for the ERRP and the discussed literature on shortest paths for EVs (\cite{Artmeier-LNAI2010,Eisner2011,Storandt2012}), the most economical route between two locations depends on the starting battery level of the vehicle (as shown in \cite{Artmeier-LNAI2010}), making routing decisions dynamic. Therefore, these problems cannot be formulated on a complete graph, since the abstraction of an arc between two retailers representing the amount of discharge (and charge) of the battery cannot be made. The graph structure of the ERRP and these works must represent the \emph{original road}, allowing for \emph{turn-by-turn} directions. 
	

	
	Another main difference of our work with respect VRP and green logistics literature is on the nature of demand. In these problems, the demand of each retailer node is fixed. In the ERRP, retailers face dynamic demand periodically, defined over a discrete finite time horizon.  Turn-by-turn decisions of the vehicle are taken over the same time dimension, which also tracks the battery level and the load inventory of the vehicle. This allows the vehicle to visit a retailer or transit a particular arc multiple times if needed, unlike in capacitated VRP problems, where visiting each retailer only once preserves the value of load inventory decision variables. The literature regarding transportation problems in which the system changes dynamically as the vehicle visits retailers is scarce. One example of such literature is the  Dependent VRP \cite{MalandrakiTDVRP}, where travelling times depends not only on distance but on time of the day. Another recent example is the Dynamic Bowser Routing Problem \cite{Rossi-DBRP}, where a bowser delivers fuel to construction machinery (assets), which can change location on the graph with each review of the system. Having discussed and motivated our approach by establishing the differences with related literature, we next proceed and describe the problem setup. 

    \begin{table}[]
    \caption{Literature overview *Recharging can happen only on downhill road segments by regenerative braking}
    \label{table:literature}
    \scalebox{0.78}{
    \begin{tabular}{llccccc}
                                                                      & Ref.                                             & \begin{tabular}[c]{@{}c@{}}Load affects \\ travel costs\end{tabular} & \begin{tabular}[c]{@{}c@{}}Recharging\\ decisions\end{tabular} & \begin{tabular}[c]{@{}c@{}}Product\\  inventory\end{tabular}     & \begin{tabular}[c]{@{}c@{}}Graph\\  structure\end{tabular} & \begin{tabular}[c]{@{}c@{}}Time \\ periods\end{tabular} \\
(C)VRP(TW)                                                            & \cite{Dantzig-truck,CVRPalgorithms,VRPTWalgorithms}              &                                                                      &                                                                &                                                                  & Complete                                                   &                                                         \\
PRP                                                                   & \cite{PRP-2011}                                  & \checkmark                                                           &                                                                & Vehicle                                                          & Complete                                                   &                                                         \\
G-VRP                                                                 & \cite{Erdogan2012}                               &                                                                      & Nodes                                                          &                                                                  & Complete                                                   &                                                         \\
E-VRTPTW                                                              & \cite{Schneider2014}                             &                                                                      & Nodes                                                          & Vehicle                                                          & Complete                                                   &                                                         \\
\begin{tabular}[c]{@{}l@{}}E-VRTPTW \\ energy model\end{tabular}          & \cite{Goeke2015}                                 & \checkmark                                                           & Nodes                                                          & Vehicle                                                          & Complete                                                   &                                                         \\
TDEVRP                                                                & \cite{Lu2020}                                    &                                                                      & Nodes                                                          & Vehicle                                                          & Complete                                                   &                                                         \\
\begin{tabular}[c]{@{}l@{}}E-VRPTW \\ partial recharging\end{tabular} & \cite{Schiffer2017}                              &                                                                      & Nodes                                                          & Vehicle                                                          & Complete                                                   &                                                         \\
\begin{tabular}[c]{@{}l@{}}E-VRP TD \\ queueing\end{tabular}          & \cite{Keskin2019}                                &                                                                      & Nodes                                                          & Vehicle                                                          & Complete                                                   &                                                         \\
TDVRP                                                                 & \cite{MalandrakiTDVRP}                           &                                                                      & Nodes                                                          & Vehicle                                                          & Complete                                                   & \checkmark                                              \\
IRPs                                                                  & \cite{IRPsReview}                                &                                                                      &                                                                & \begin{tabular}[c]{@{}c@{}} Vehicle \& \\  Retailer\end{tabular} & Complete                                                   & \checkmark                                              \\
\begin{tabular}[c]{@{}l@{}}EVs shortest\\  paths\end{tabular}         & \cite{Artmeier-LNAI2010,Eisner2011,Storandt2012} &                                                                      & Arcs*                                                          &                                                                  & original-road                                              &                                                         \\
DBRP                                                                  & \cite{Rossi-DBRP}                                &                                                                      &                                                                & \begin{tabular}[c]{@{}c@{}} Vehicle \& \\  Retailer\end{tabular} & original-road                                              & \checkmark                                              \\
ERRP                                                                  &                                                  & \checkmark                                                           & Arcs                                                           & \begin{tabular}[c]{@{}c@{}} Vehicle \& \\  Retailer\end{tabular} & original-road                                              & \checkmark                                             
\end{tabular}}
    \end{table}

	\newpage
	
	\section{Problem definition and formulation}
	\label{sec:ERRP}

	We present the Electric Roads Routing Problem (ERRP) as a transportation problem in which some road sections feature ERS technology, allowing a single vehicle to be powered by grid electricity while charging its electric battery alongside. The aim of the vehicle is to deliver a single product to a set of retailers facing demand of that product. In this section we first provide the problem definition, and then introduce a SDP formulation.
	
	\subsection{Problem definition}
	\label{sec:ERRPdef}
	We consider a graph $G = \langle \mathcal{N},\mathcal{A} \rangle$ representing a road network, with $\mathcal{N}$ as the set of nodes and $\mathcal{A}$ the set of arcs as segments of road.
	\noindent
    A subset of the nodes  $\mathcal{C}\subseteq \mathcal{N}$ represents retailers, and there is a node, denoted as 0, representing a depot.
    \noindent
    Other nodes in $\mathcal{N} \setminus \mathcal{C}\cup \{0\}$ represent road intersections or are used to delimit ERS road sections.

    We consider a discrete planning horizon of $T$ periods. 
    \noindent
    At each period $t$, each retailer $c\in \mathcal{C}$ faces a demand of a single product $d^c_t$.
    \noindent
    A vehicle movement from its current node position to an adjacent one takes one time period.
    \noindent
    If the origin node is the depot or a retailer, the vehicle can load or deliver some inventory, respectively. Each retailer $c\in \mathcal{C}$ has an inventory capacity of $k_c$ units, and the vehicle has a hold capacity for inventory of $K$ units.
    \noindent
    For each period, after a possible delivery to a retailer, the demand of all retailers takes place. If a retailer does not hold enough inventory to cover the demand of that period, the unsatisfied demand is considered lost and a penalty cost $p$ per unit is incurred. 

    The vehicle has an electric battery with an energy capacity of $B$ (in kWh units). 
    \noindent
    When the vehicle traverses an arc $(i,j)$ it requires $r_{ij}(M)= \alpha_{ij}M+\beta_{ij}$ kWh from the battery, where $M$ is unladen weight of the vehicle (in kg.) and $\alpha_{ij}$ and $\beta_{ij}$ are known parameters (see Appendix \ref{sec:energymodel}).
    \noindent
    Parameters $s_{ij}$ account for the energy supplied  to the vehicle (in kWh) when it traverses an electrified arc. 
    \noindent
    The level of electric battery of the vehicle is updated after traversing an arc, accounting for the supplied and required electric energy at the time. If the level of battery is not sufficient the vehicle uses fuel (e.g. gasoline) to reach the destination node; the quantity of fuel consumed is computed on the basis of the extra energy required to arrive to the new node, after the battery is depleted.
    \noindent
    Although recharging opportunities present only on ERS arcs, it would be straightforward to extend the model to include charging stations: spending one or more periods at a charging station node would charge the battery according to the power the station can provide and the granularity of the periods.
    \noindent
    The cost of a kWh of energy provided by grid electricity or battery and fuel is given by parameters $C^e$ and $C^f$ respectively, and they are used to derive the transportation costs of the vehicle.
    
    The goal of the ERRP is to find a route for the vehicle, and a replenishment plan for retailers, which minimises the costs of travelling, and the costs of unfulfilled demand over the time horizon. 

    
    \begin{figure}[ht]
    \caption{Problem instance on a ERS network. Tuples on arcs are $\langle \alpha_{ij}, \beta_{ij} \rangle$; Supplied energy in the dashed blue arc is $s_{04}=20$, while the rest are zero.}
    \label{fig:graph_example}
    \centering
    \begin{tikzpicture}[scale=1,transform shape]
      \Vertex[x=0,y=0]{0}
      \Vertex[x=4.5,y=2.6]{1}
      \Vertex[x=4.5,y=-2.6]{2}
      \Vertex[x=2.25,y=1.3]{3}
      \Vertex[x=3,y=0]{4}
      \Vertex[x=2.25,y=-1.3]{5}
      \tikzstyle{LabelStyle}=[fill=white,sloped]
      \Edge[label={$\langle1,1\rangle$}](0)(3)
      \Edge[label={$\langle2,1\rangle$},style={color=cblue,dashed},](0)(4)
      \Edge[label={$\langle1,1\rangle$}](0)(5)
      \Edge[label={$\langle2,1\rangle$}](4)(1)
      \Edge[label={$\langle2,1\rangle$}](4)(2)
      \Edge[label={$\langle1,1\rangle$}](3)(1)
      \Edge[label={$\langle1,1\rangle$}](1)(2)
      \Edge[label={$\langle1,1\rangle$}](2)(5)
    \end{tikzpicture}
    \end{figure}
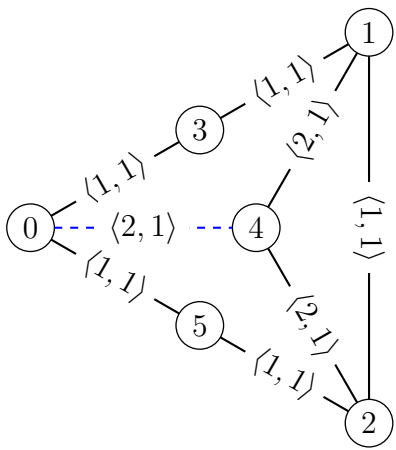

	
	\begin{example}
	\label{sec:example}
	The following example is presented with the illustrative purposes of depicting the ERRP. For the sake of simplicity, the values of the parameters are abstract.

    \begin{table}[h]
        \caption{Optimal solution for Example \ref{sec:example}}
        \label{table:exampleSol1}
        \begin{tabular}{llllll}
                            & t = 0 & t = 1 & t = 2 & t = 3 & t = 4 \\
        V. Position         & N/A   & 0     & 4     & 1     & 2     \\
        V. inv. upload      & N/A   & 3     & 0     & 0     & 0     \\
        Delivery            & N/A   & 0     & 0     & 2     & 1     \\
        Battery level       & N/A   & 0     & 11    & 2     & 0     \\
        Weight $M = w + L$  & N/A   & 4     & 4     & 2     & 1     \\
        Vehicle Load (L)    & 0     & 3     & 3     & 1     & 0     \\
        Inv. ret. 1         & 2     & 1     & 0     & 1     & 0     \\
        Inv. ret. 2         & 3     & 2     & 1     & 0     & 0     \\
        Required energy     & N/A   & 9     & 9     & 3     & 0     \\
        Travel costs        & N/A   & 9     & 9     & 7     & 0     \\
        Penalty costs       & N/A   & 0     & 0     & 0     & 0    
        \end{tabular}
    \end{table}
	
    We consider the simple road network presented in Figure \ref{fig:graph_example} where the depot is located at node 0 and the only retailers of the network are located at nodes 1 and 2, with inventory capacities $k_1 = k_2 = 5$. Demand is $d^c_t=1$ for $c =1,2$ and $t=1,\ldots,T$ and the penalty cost per unit of lost sales is $p=25$. The unladen weight of the vehicle is $w=1$, and has a battery capacity of $B = 20$ and inventory capacity $K = 4$. 
    Travel costs are $C^e = 1$ and $C^f = 5$. 
    We consider a time horizon of $T= 4$ periods. The vehicle starts at the depot node in period 1, with no load and no battery. Retailer 1 starts with 3 units of inventory and Retailer 2 with 2 units.

	The optimal solution for the problem is presented in Table \ref{table:exampleSol1}. On the first period, the vehicle starts at the depot and loads 3 units of inventory. Then transits to node 3 using a ERS road. In period 2, the vehicle starts with 9 units of battery (on the previous travel the ERS road provided 20 units but 9 were required to power the vehicle) and transits to node 1. In period 3, the vehicle delivers 2 units of inventory to Retailer 1 and transits to node 2 using 2 units of battery and 1 of fuel. In period 4 the vehicle delivers inventory to Retailer 2, and the vehicle does not transit to a new node. There are no penalty costs incurred in any periods since the vehicle delivers to both retailers before demands exceeds inventory and enough to cover demand until the end of the time horizon. 

    The optimal solution when $(0,4)$ is not electrified follows similarly and is presented in Table \ref{table:exampleSol2}. In this case the path of the vehicle does not transit $(0,4)$ since the arc is not electrified and travelling costs following the former path would be more expensive. 
	
    \begin{table}[h]
        \caption{Optimal solution for Example \ref{sec:example} with $s_{04}=0$}
        \label{table:exampleSol2}
        \begin{tabular}{llllll}
                          & t = 0 & t = 1 & t = 2 & t = 3 & t = 4 \\
        V. Position             & N/A   & 0     & 3     & 1     & 2     \\
        V. inv. upload          & N/A   & 3     & 0     & 0     & 0     \\
        Delivery                & N/A   & 0     & 0     & 2     & 1     \\
        Battery level           & N/A   & 0     & 0     & 0     & 0     \\
        Weight $M = w + L$      & N/A   & 4     & 4     & 2     & 1     \\
        Vehicle Load (L)        & 0     & 3     & 3     & 1     & 0     \\
        Inv. Ret. 1             & 2     & 1     & 0     & 1     & 0     \\
        Inv. Ret. 2             & 3     & 2     & 1     & 0     & 0     \\
        Required energy         & N/A   & 5     & 5     & 3     & 0     \\
        Travel costs            & N/A   & 25    & 25    & 18    & 0     \\
        Penalty costs           & N/A   & 0     & 0     & 0     & 0    
        \end{tabular}
    \end{table}
	
	\end{example}
	
	\section{A Stochastic Dynamic Programming formulation}
    \label{sec:SDP-SERRP}
    
     In what follows, we will consider a stochastic variant of the ERRP, in which  demand is stochastic. Demand uncertainty affects vehicle load and delivery decisions, and therefore the total weight of the vehicle in each period, which is one of the main factors affecting energy requirements as discussed in \cite{PRP-2011}. We refer to this variant of the problem as the Stochastic ERRP (S-ERRP). The  S-ERRP can be formulated via stochastic dynamic programming (SDP) over a planning horizon of $T$ periods by defining the \emph{state space} of the problem, the \emph{feasible actions} $A_s$ associated to each state $s$, the \emph{transition probabilities} $p^a_{s,s'}$ and \emph{immediate costs} $c^a_s$ of each state-action pair $\langle s,a \rangle$ and the \emph{objective function} defined over the value function:
    
    \begin{equation}
        \label{eq:value-function_Bellmans}
        V_t(s) = \min_{k \in A_s} \{c^k_s + \sum_{s'} p^k_{s,s'} V_{t+1}(s')  \}
    \end{equation}

    Before defining these components we provide a detailed explanation of the dynamics of the vehicle and the retailers.
	\noindent
	At the beginning of each period the vehicle can deliver product to a retailer, if the starting  node is a retailer node, or load up inventory to the vehicle, if the starting node is the depot node. The total weight of the vehicle for the period is updated after these possible actions.
	\noindent
	In both cases the quantities are subject to the storage capacity of the retailers $k_c$ and the vehicle inventory capacity $K$.
	\noindent
	After any possible delivery during the period is accounted for, the demand on each retailer (modeled by random variables $d^i_t$ for retailer $i$ and period $t$) takes place. 
	\noindent
	If a retailer does not hold enough inventory to cover the demand of each period, the unsatisfied demand of that period is considered lost and a penalty cost $p$ per unit is incurred. 
    \noindent
    The vehicle finally transits to a new node reaching an adjacent node to the origin one by the end of the period, updating the level of battery for the next period and incurring travelling costs.
    \noindent
    We next model the S-ERRP by using SDP:
    
    \begin{itemize}
        \item \emph{States}: States are given by tuples $(t,V_{\text{pos}},V_{\text{inv}},V_{\text{bat}},I_{\text{ret}})$ where $t$ is the current period associated with the state, $V_{\text{pos}}$ is the starting position of the vehicle at time $t$, $V_{\text{inv}}$ is the starting product inventory on hold of the vehicle before any load or delivery of product takes place, $V_{\text{bat}}$ the battery level of the vehicle at $V_{\text{pos}}$ and $I_{\text{ret}}$ is a tuple with the inventory of each retailer before any delivery at time $t$.
        \item \emph{Actions}: An action $a\in A_s$, related to state $s$, is a tuple $(V'_{\text{pos}},V_{\text{up}},Q_{\text{ret}})$ where $V'_{\text{pos}}$ is the destination node of the vehicle, $V_{\text{up}}$ is the inventory loaded up to the vehicle and $Q_{\text{ret}}$ is the possible delivery quantity at $V_{\text{pos}}$.
        \item \emph{Transition probabilities}: From a state $s$ and action $a$, the transition probabilities to a new state $s'$, $p^a_{s,s'}$ are derived directly from the random demand distributions.
        \item \emph{Immediate costs}: The immediate costs of a state-action pair $(s,a)$ are given by the costs of travel from $V_{\text{pos}}$ to $V'_{\text{pos}}$ with weight derived from $V_{\text{inv}}$ and accounting for the weight that is loaded up or delivered to a retailer, and the expected penalty costs after updating $I_{\text{ret}}$ with the possible delivery $Q_{\text{ret}}$ at node $V_{\text{pos}}$.
    \end{itemize}
    
    With all components defined, the optimal policy for each state can be derived by applying a backward of forward recursion algorithm on the value function (Eq. \ref{eq:value-function_Bellmans}). However, this approach only allow us to solve small instances of the problem. In the next section, we explore a MILP heuristic for the S-ERRP.

    \section{A stochastic MILP heuristic}
    \label{sec:MILPformulation}

    In this section, we present a mixed-integer linear programming (MILP) heuristic for the S-ERRP. A summary of the symbols used in this model is presented in Appendix \ref{sec:MILPsymbols}. We define the adjacency matrix of graph $\mathcal{G}$ with the binary parameters $\delta_{ij}$. We also denote $N$ and $C$ as the cardinals of $\mathcal{N}$ and $\mathcal{C}$, the number of nodes and retailers.

    In the S-ERRP, the demand faced by each retailer $i$ and period $t$ is a random variable denoted as $d^i_t$. As a result, the lost sales of each retailer and period are stochastic variables that depend on the deliveries given, and become a central aspect when formulating stochastic programming models in inventory control. Given a random variable $\omega$ and a scalar $q$, the first order loss function is defined as
    
    \begin{equation}
        \label{eq:lossfunction}
        \mathcal{L}_\omega(q) = E[\max (\omega - q, 0)]
    \end{equation}
    In a stochastic lot sizing setting, it expresses the expected value of the exceeding demand, given starting inventory at hand $q$, when demand follows $\omega$. Reciprocally, the complementary first order loss function expresses the expected level of inventory after demand:
    
    \begin{equation}
        \label{eq:complossfunction}
        \mathcal{\hat{L}}_\omega(q) = E[\max (q - \omega, 0)]
    \end{equation}
    
    In what follows, we follow a policy considering that replenishment periods and replenishment quantity decisions for all retailers are fixed at the beginning of the planning horizon, before demand random variables for any period are realised. In practical terms this means that the route the vehicle follows is decided beforehand as well, and since the replenishment quantities $Q^i_t$ are fixed, the load of the vehicle and energy requirements can be calculated at the beginning of the planning horizon. 
     
    The objective function minimises the total carbon emissions emitted during the time horizon and the cost of lost sales at the retailers:
	
    \begin{equation}
    \label{eq:stochastic_obj}
        \text{min} \sum_{t =2}^T \sum_{i=1}^N \sum_{j=1}^N C^e s_{ij} T^{ij}_{t-1} + \sum_{t =2}^T C^e E^b_t + C^f E^f_t + \sum_{t=1}^T \sum_{i=1}^C p[I^i_t]^-
    \end{equation}
    

    The first term of the objective function accounts for the emissions caused by powering the vehicle on transited ERS roads. 
    \noindent
    The second and third term account for the emissions related to the used battery and the used fuel respectively.
    \noindent
    The last term accounts for the penalty cost of the expected lost sales.

	
	We follow introducing constraints related to the movement of the vehicle, product inventory and energy use of the vehicle. The first set of Eqs. (\ref{eq:electric_one_position}-\ref{eq:electric_TVlink4}) deal with the movement of the vehicle over the time horizon.
	\noindent
	At each period of the planning horizon, the vehicle is situated at one and only one of the nodes $i \in \mathcal{N}$. We model this by considering binary variables $V^i_t$, being equal to one if the vehicle is present at node $i \in \mathcal{N}$ at the beginning of period $t$, and $T^{ij}_t$ as a binary variable set to one if the vehicle transits from node $i$ to node $j$ during period $t$. 
	\noindent
	Eqs. \ref{eq:electric_one_position} reflects the fact that the vehicle is, at any point in time, at one and only one node of graph $\mathcal{G}$.
	\noindent
	Eqs. \ref{eq:electric_transit} ensures that for each period the vehicle transits to a new node from where it starts.
	\noindent
	Eqs. \ref{eq:electric_TVlink1}-\ref{eq:electric_TVlink3} link transit $T^{ij}_t$ and position $V^i_t$ variables.
	\noindent
	Eqs. \ref{eq:electric_TVlink4} ensures that the vehicle only transits valid arcs, observing the adjacency matrix of the graph.
	\noindent

    \begin{align}
        \label{eq:electric_one_position}
        \sum_{i=1}^N V^i_t = 1 & &  t = 1,\ldots,T\\
        \label{eq:electric_transit}
        \sum_{j=1}^N T^{ij}_{t-1} = V^i_{t-1} & & i = 0,\ldots,N-1; t = 2,\ldots,T\\
        \label{eq:electric_TVlink1}
        T^{ij}_{t-1} \geq V^i_{t-1} + V^j_t - 1 & &  i,j = i = 0,\ldots,N-1; t = 2,\ldots,T\\
        \label{eq:electric_TVlink2}
        T^{ij}_{t-1} \leq V^i_{t-1} & &  i,j = i = 0,\ldots,N-1; t = 2,\ldots,T\\
        \label{eq:electric_TVlink3}
        T^{ij}_{t-1} \leq V^j_t & &  i,j = i = 0,\ldots,N-1; t = 2,\ldots,T\\
        \label{eq:electric_TVlink4}
        T^{ij}_{t} \leq \delta_{ij} & & i,j = 0,\ldots,N-1; t = 1,\ldots,T-1  
    \end{align}

    The second set of Eqs. (\ref{eq:filling_truck}-\ref{eq:weight_vehicle}) are regarding with product inventory levels: Eqs. \ref{eq:filling_truck} ensure the vehicle can only load inventory to the vehicle while visiting the depot node, and only up to the vehicle capacity $K$. Eqs. \ref{eq:delivered_demand} ensure the vehicle can only deliver inventory in one period if it is visiting a retailer node, and up to their capacity. The transition of inventory levels of the vehicle are tracked by Eqs. \ref{eq:cargo1} and \ref{eq:cargo2}, ensuring the total cargo does not exceed the capacity of the vehicle. 
    \noindent
    Eqs. \ref{eq:expected_shortage} approximate expected lost sales $[I^i_t]^-$ for each retailer and period using the loss function using the random variables $\mathbf{d}^i_{1t} = \mathbf{d}^i_1+\ldots+\mathbf{d}^i_t$. The idea of the approximation is to calculate the lost sales of a single period comprising the demand of periods $1,\ldots,t$ with an starting inventory of the starting inventory $s_i$, and all deliveries up to period $t$, plus the losses until period $t-1$ (see Lemma 5 from \cite{Rossi-DBRP}), but also accounting for a reduction of the starting inventory of the expected exceeding capacity deliveries of periods $1,\ldots,t$.
    \noindent
    Similarly, Eqs. \ref{eq:expected_inventory} approximate expected inventory of each retailer and period $[I^i_t]^+$.
    \noindent
    The loss function and the complementary loss function are not linear. We follow the approach of \cite{Rossi-linearloss} to linearise these functions to be able to use them in our stochastic programming. model. 
    \noindent
    Eqs. \ref{eq:exceeding_capacity} are used to calculate the expected inventory that would exceed the capacity of each retailer when a delivery is made. It follows from the fact that in a lost sales setting shortages are reset at the end of each period and $Q^i_t$ and $s_i$ are constants (Lemma 6 from \cite{Rossi-DBRP}).
    \noindent
    Finally, Eqs. \ref{eq:weight_vehicle} update the total weight of the vehicle at the end of each period, after a possible delivery or load up of inventory.

    \begin{align}
        \label{eq:filling_truck}
        L_t \leq K V^0_t & & t = 1,\ldots,T\\
        \label{eq:delivered_demand}
        Q^i_t \leq k_i  V^i_t & &  i = 1,\ldots,C; t = 1,\ldots,T\\
        \label{eq:cargo1}
        l + \sum_{k=1}^t L_k - \sum_{i = 1}^N \sum_{k=1}^{t-1} Q^i_k \leq K  & & t = 1,\ldots,T \\
        \label{eq:cargo2}
        l + \sum_{k=1}^t L_k - \sum_{i = 1}^N \sum_{k=1}^{t-1} Q^i_k \geq 0  & & t = 1,\ldots,T \\
        \label{eq:expected_shortage}
        [I^i_t]^- = \mathcal{L}_{d^i_{1t}}\left(s_i + \sum_{k=1}^t Q^i_k + \sum_{k=1}^{t-1}[I^i_{t-1}]^- - \sum_{k=1}^t [E^i_k] \right) & & t = 1,\ldots,T; i = 1,\ldots, C  \\
        \label{eq:expected_inventory}
        [I^i_t]^+ = \mathcal{\hat{L}}_{d^i_{1t}}\left(s_i + \sum_{k=1}^t Q^i_k + \sum_{k=1}^{t-1}[I^i_{t-1}]^- - \sum_{k=1}^t [E^i_k] \right) & & t = 1,\ldots,T; i = 1,\ldots, C  \\
        \label{eq:exceeding_capacity}
        [E^i_t] = \max \left( [I^i_t]^+ + Q^i_t - s_i,0\right) & & t = 1,\ldots,T; i = 1,\ldots, C \\  
        \label{eq:weight_vehicle}
        W_t = w + l + \sum_{k=1}^t L_k - \sum_{i = 1}^N \sum_{k=1}^{t} Q^i_k  & & t = 1,\ldots,T
    \end{align}

    The set of Eqs. (\ref{eq:unbounded_flow}-\ref{fuel_energy}) track the battery level inventory over the time horizon, which is used to derive the energy used from battery and fuel for each period. 
    \noindent
    The variables $E^b_t$ represent the energy used by the vehicle (in kWh) from its battery on the travel performed during period $t$. In contrast, variables $E^f_t$ represent the mechanical energy the vehicle needs to finish traversing an arc when the battery has been depleted.
	\noindent
	In order to determine these variables, we define further auxiliary variables and parameters.
    \noindent
    First of all, variables $b^u_t$ represent the unconstrained battery level of the vehicle, calculated in Eqs. \ref{eq:unbounded_flow} from the actual level of battery of the previous period, $b_{t-1}$, plus and minus supplied and required energy from the total mass of the vehicle $W_{t-1}$ at the beginning of the period. Expanding this equation, the product of a binary and a continuous variable appear, $T^{ij}_{t-1}W_{t-1}$, which are linearised in our model considering $L_t$ variables, the weight of the vehicle at the end of period $t$, are bounded between zero and the maximum weight of the vehicle $W$. In Eqs. \ref{eq:battery_level}, variables $b_t$ are calculated as the $b^u_t$ bounded between 0 and $B$, representing the actual level of battery for each period.
    \noindent
	Eqs. \ref{eq:battery_usage} calculate battery usage in each period, considering that if the battery was charged traversing a ERS arc, the battery consumed during that period should be zero.
	\noindent
	In Eqs. \ref{fuel_energy}, if $b^u_t<0$, the interpretation is that during period $t$ the battery was depleted and $-b^u_t$ kWh were still needed to complete travelling the arc on battery. Dividing by the corresponding battery energy losses factor $\lambda$, we obtain $E^f_t$ variables as the required mechanical energy to complete the travel when battery is not sufficient.

    \begin{align}
        \label{eq:unbounded_flow}
        b^u_t = b_{t-1} + \sum_{i=1}^N \sum_{j=1}^N T^{ij}_{t-1} \left(s_{ij} - \alpha_{ij} W_{t-1} - \beta_{ij}\right) && t = 2,\ldots,T\\
        \label{eq:battery_level}
        b_t = \text{min}\{\text{max}\{b^u_t,0\},B \} && t = 2,\ldots,T\\
        \label{eq:battery_usage}
        E^b_{t-1} = \text{max}\{b_{t-1} - b_t,0\} && t = 2,\ldots,T\\
        \label{fuel_energy}
        E^f_{t-1} = \max\left(-b^u_t,0\right)/\lambda && t = 2,\ldots,T
    \end{align}

    Finally, Eqs. \ref{eq:domains}-\ref{eq:lin_domains4} represent the domains of the variables used.

    \begin{align}	
        \label{eq:domains}
        T^{ij}_t, V^i_t \in \{0,1\} & & c = 1,\ldots, C; t = 1,\ldots,T\\
        \label{eq:lin_domains1}
	    b_t,b^u_t,E^f_t,E^b_t, L_t \geq 0 & & t = 2,\ldots,T\\
	    \label{eq:lin_domains2}
	    b_t \leq B & & t = 2,\ldots,T\\
	    \label{eq:lin_domains3}
	    Q^i_t, S^i_t \geq 0 & & i = 1,\ldots, C; t = 1,\ldots,T\\
	    \label{eq:lin_domains4}
	    0 \leq S^i_t \leq d^i_t & & i = 1,\ldots, C; t = 1,\ldots,T
	\end{align}

The MILP model consists of the objective function \ref{eq:stochastic_obj} subject to the constraints set in Eqs. (\ref{eq:electric_one_position}-\ref{eq:lin_domains4})

    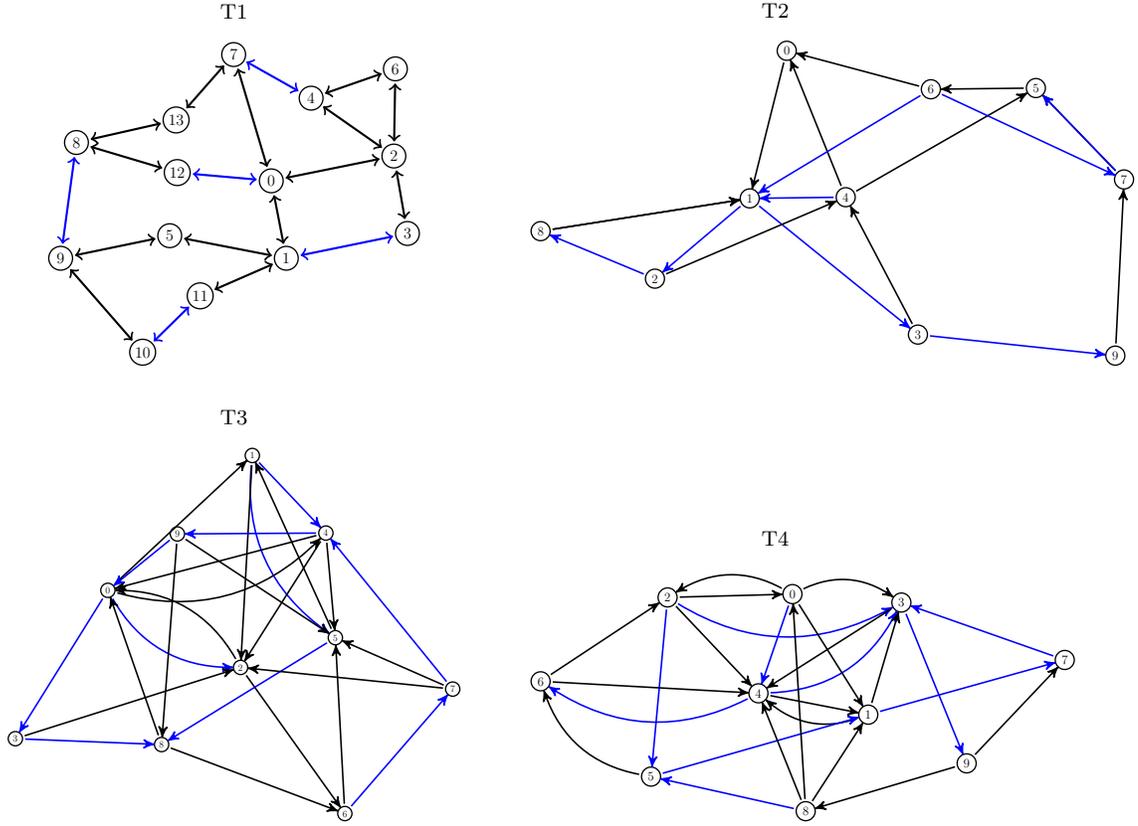
\begin{figure*}[ht]
    \caption{Network topologies}
    \captionsetup[subfigure]{labelformat=empty}
    \label{fig:computational_exp_graph}
    \centering
    \begin{subfigure}[b]{0.475\textwidth}
        \centering
        \caption{T1}
        \begin{tikzpicture}[scale=0.5, transform shape]
          \Vertex[x=11.2,y=7.24]{0}
          \Vertex[x=11.6,y=5.18]{1}
          \Vertex[x=14.46,y=7.9]{2}
          \Vertex[x=14.82,y=5.84]{3}
          \Vertex[x=12.26,y=9.44]{4}
          \Vertex[x=8.5,y=5.78]{5}
          \Vertex[x=14.5,y=10.22]{6}
          \Vertex[x=10.2,y=10.6]{7}
          \Vertex[x=6.02,y=8.26]{8}
          \Vertex[x=5.6,y=5.18]{9}
          \Vertex[x=7.78,y=2.68]{10}
          \Vertex[x=9.31,y=4.18]{11}
          \Vertex[x=8.7,y=7.46]{12}
          \Vertex[x=8.67,y=8.86]{13}
          \tikzstyle{EdgeStyle}=[pre and post]
          \Edges[](0, 1)
          \Edges[style={color=cblue,}](1, 3)
          \Edges[](1, 11)
          \Edges[](2, 3)
          \Edges[](4, 6)
          \Edges[](7, 13) 
          \Edges[](8, 13)
          \Edges[style={color=cblue,}](10, 11)  
          \Edges[](0, 2)
          \Edges[](0, 7)
          \Edges[style={color=cblue,}](0, 12)
          \Edges[](1, 5)
          \Edges[](2, 4)
          \Edges[](2, 6)
          \Edges[style={color=cblue,}](4, 7)
          \Edges[](5, 9)
          \Edges[style={color=cblue,}](8, 9)
          \Edges[](8, 12)
          \Edges[](9, 10)          
        \end{tikzpicture}
    \end{subfigure}
    \hfill
        \begin{subfigure}[b]{0.475\textwidth}
        \centering
        \caption{T2}
        \begin{tikzpicture}[scale=0.4, transform shape]
          \Vertex[x=9.78,y=12.67]{0}
          \Vertex[x=8.55,y=7.76]{1}
          \Vertex[x=5.4,y=5.09]{2}
          \Vertex[x=14.14,y=3.23]{3}
          \Vertex[x=11.74,y=7.8]{4}
          \Vertex[x=18.06,y=11.43]{5}
          \Vertex[x=14.57,y=11.39]{6}
          \Vertex[x=20.99,y=8.39]{7}
          \Vertex[x=1.6,y=6.67]{8}
          \Vertex[x=20.7,y=2.53]{9}
          \tikzstyle{EdgeStyle}=[post]
          \Edges[](0, 1)
          \Edges[style={color=cblue,}](1, 2)
          \Edges[style={color=cblue,}](1, 3)
          \Edges[](2, 4)
          \Edges[style={color=cblue,}](2, 8)
          \Edges[](3, 4)
          \Edges[style={color=cblue,}](3, 9)
          \Edges[](4, 0)
          \Edges[style={color=cblue,}](4, 1)
          \Edges[](4, 5)
          \Edges[](5, 6)
          \Edges[](6, 0)
          \Edges[style={color=cblue,}](6, 1)
          \Edges[style={color=cblue,}](6, 7)
          \Edges[](7, 5)
          \Edges[](8, 1)
          \Edges[style={color=cblue,}](7, 5)  
          \Edges[](8, 1)
          \Edges[](9, 7)
        \end{tikzpicture}
    \end{subfigure}
    \vskip\baselineskip
    \begin{subfigure}[b]{0.475\textwidth}    
        \centering
        \caption{T3}
        \begin{tikzpicture}[scale=0.3, transform shape]
          \Vertex[x=5.34,y=10.46]{0}
          \Vertex[x=11.74,y=16.44]{1}
          \Vertex[x=11.22,y=7.04]{2}
          \Vertex[x=1.24,y=3.88]{3}
          \Vertex[x=15,y=13]{4}
          \Vertex[x=15.42,y=8.36]{5}
          \Vertex[x=15.84,y=0.56]{6}
          \Vertex[x=20.62,y=6.08]{7}
          \Vertex[x=7.72,y=3.62]{8}
          \Vertex[x=8.42,y=12.96]{9}
          \tikzstyle{EdgeStyle}=[post,bend right]
          \Edges[style={color=cblue,}](1, 5)
          \Edges[style={color=cblue,}](0, 2)
          \Edges[](0, 4)
          \Edges[](2, 0)
          \tikzstyle{EdgeStyle}=[post]
          \Edges[](0, 1)
          \Edges[style={color=cblue,}](0, 3)
          \Edges[](1, 2)
          \Edges[style={color=cblue,}](1, 4)
          \Edges[](2, 6)
          \Edges[](3, 2)
          \Edges[style={color=cblue,}](3, 8)
          \Edges[](4, 0)
          \Edges[](4, 2)
          \Edges[](4, 5)
          \Edges[style={color=cblue,}](4, 9)
          \Edges[](5, 1)
          \Edges[style={color=cblue,}](5, 8)
          \Edges[](6, 5)
          \Edges[style={color=cblue,}](6, 7)
          \Edges[](7, 2)
          \Edges[style={color=cblue,}](7, 4)
          \Edges[](7, 5)
          \Edges[](8, 0)
          \Edges[](8, 6)
          \Edges[style={color=cblue,}](9, 0)
          \Edges[](9, 5)
          \Edges[](9, 8)
        \end{tikzpicture}
    \end{subfigure}    
    \hfill
    \begin{subfigure}[b]{0.475\textwidth}
        \centering
        \caption{T4}
        \begin{tikzpicture}[scale=0.40, transform shape]
          \Vertex[x=10.56,y=9.6]{0}
          \Vertex[x=13.08,y=5.59]{1}
          \Vertex[x=6.41,y=9.49]{2}
          \Vertex[x=14.18,y=9.32]{3}
          \Vertex[x=9.43,y=6.3]{4}
          \Vertex[x=5.86,y=3.53]{5}
          \Vertex[x=2.19,y=6.7]{6}
          \Vertex[x=19.61,y=7.41]{7}
          \Vertex[x=11,y=2.38]{8}
          \Vertex[x=16.35,y=3.97]{9}  
          \tikzstyle{EdgeStyle}=[post,bend right]
          \Edges[](0, 2)
          \Edges[style={color=cblue,}](2, 3)
          \Edges[style={color=cblue,}](4, 3)
          \tikzstyle{EdgeStyle}=[post,bend left]
          \Edges[](0, 3)
          \Edges[](1, 4)
          \Edges[style={color=cblue,}](4, 6)
          \Edges[](5, 6)
          \tikzstyle{EdgeStyle}=[post]
          \Edges[](0, 1)
          \Edges[style={color=cblue,}](0, 4)
          \Edges[](1, 3)
          \Edges[style={color=cblue,}](1, 7)
          \Edges[](2, 0)
          \Edges[](2, 4)
          \Edges[style={color=cblue,}](2, 5)
          \Edges[](3, 4)
          \Edges[style={color=cblue,}](3, 9)
          \Edges[](4, 1)
          \Edges[style={color=cblue,}](5, 1)
          \Edges[](6, 2)
          \Edges[](6, 4)
          \Edges[style={color=cblue,}](7, 3)
          \Edges[](8, 0)
          \Edges[](8, 1)
          \Edges[](8, 4)
          \Edges[style={color=cblue,}](8, 5)
          \Edges[](9, 7)
          \Edges[](9, 8)
        \end{tikzpicture}
    \end{subfigure}
    \end{figure*}

    \section{Numerical experiments: effectiveness of the MILP heuristic}
    \label{sec:comp_study}
    
    With the aim of investigating the effectiveness of the MILP heuristic designed for the ERRP we present the following computational study. Our goal is to compare the cost of the optimal plan obtained by the heuristic against the optimal cost obtained via SDP. 

    To ensure that the state space of the SDP implementation is finite, we introduce the following assumptions in relation to inventory and battery state variables. Please note that this is done without loss of generality, and for the sake of comparability, as these assumptions are not required in the MILP model.
    \noindent
    First, the stochastic demand faced by retailers is expressed using finite discrete random variables. Both retailers and vehicle product inventory are discretised by using the same units.
    \noindent
    Second, the battery of the vehicle needs to be also discretised, and as a consequence, the required energy will be expressed in terms of these battery levels. Equation \ref{eq:battery_energy} is used to round values to the closest battery level unit, but this constraint cannot be captured in a MILP model. We will therefore rely on Constraint Programming \cite{CP-Handbook} to model and solve these instances.
    \noindent
    
    All our experiments have been carried out on a PC  Intel(R) Core(TM) i7-8665U CPU @ 1.90GHz with 32 GB of RAM, using IBM ILOG CPLEX Studio 12.7 with a time limit of 10 minutes and default constraint programming settings otherwise.
    \noindent
    The design of the experiments is described in Section \ref{sec:exp_design}, and the results obtained are presented and discussed in Section \ref{sec:results}.

    \begin{table}[h]
    \caption{Values of parameters considered in the test bed}

    \begin{tabular}{lll}
    \label{tab:testbedval}
    \begin{tabular}[c]{@{}l@{}}Initial inventory \\ at \{R1,R2\}\end{tabular} & \{\{0, 0\},\{5, 5\}\}                                                                                                                                                                                                                                                                                               &  \\
    \begin{tabular}[c]{@{}l@{}}Demand \\ distributions\end{tabular}           & \begin{tabular}[c]{@{}l@{}}(D1) $\lambda_{R1}$ = \{2,2,2,2,2,2,2,2,2\}, \\     $\lambda_{R2}$ = \{2,2,2,2,2,2,2,2,2\}\\ (D2) $\lambda_{R1}$ = \{1,1,2,2,3,3,4,4,5\}, \\   $\lambda_{R2}$ = \{5,4,4,3,3,2,2,1,1\}\\ (D3) $\lambda_{R1}$ = \{1, 1, 2, 1, 1, 2, 2, 3, 1\}, \\   $\lambda_{R2}$ = \{1, 1, 2, 1, 1, 2, 2, 3, 1\}\end{tabular} &  \\
    \begin{tabular}[c]{@{}l@{}}Unit penalty \\ cost\end{tabular}              & $p$ = \{10, 20, 30\}                                                                                                                                                                                                                                                                                                & 
    \end{tabular}
    \end{table}

        \begin{table}[]
    \caption{Pivot table of mean, median and standard deviation of percentage error (MPE, MdPE, SD respectively) of the solutions obtained by the MILP heuristic for the computational study}
    \label{table:results}
    \scalebox{0.9}{
    \begin{tabular}{llll}
                      & MPE              & MdPE    & SD      \\
    Network           &                  &         &         \\
    T1                & 3.47\%          & 3.42\% & 2.14\% \\
    T2                & 3.48\%          & 3.19\% & 2.05\% \\
    T3                & 3.18\%          & 2.53\% & 2.57\% \\
    T4                & 4.73\%          & 3.88\% & 4.12\% \\
                      &                  &         &         \\
    \multicolumn{2}{l}{Initial inv.}     &         &         \\
    (0,0)             & 3.94\%          & 3.41\% & 2.87\% \\
    (5,5)             & 3.49\%          & 3.03\% & 2.91\% \\
                      &                  &         &         \\
    Penalty           &                  &         &         \\
    10                & 1.67\%          & 1.30\% & 1.57\% \\
    20                & 3.98\%          & 3.51\% & 2.45\% \\
    30                & 5.50\%          & 5.22\% & 3.06\% \\
                      &                  &         &         \\
    \multicolumn{2}{l}{Demand   pattern} &         &         \\
    D1                & 3.70\%          & 3.39\% & 2.67\% \\
    D2                & 3.50\%          & 3.02\% & 2.71\% \\
    D3                & 3.95\%          & 3.29\% & 3.28\% \\
                      &                  &         &         \\
    General           & 3.72\%          & 3.25\% & 2.90\%
    \end{tabular}}
    \end{table}
    \subsection{Design of experiments}
    \label{sec:exp_design}

    The experiments are designed over a time horizon comprising $T=9$ periods, on four road networks depicted as the graphs in Figure \ref{fig:computational_exp_graph}. Road network T1 has been generated for this study, while the other three were used in the computational study of \cite{Rossi-DBRP}.
    \noindent
    On each road network, we consider six sets of retailers chosen randomly.
    \noindent
    The demand at the retailers follow normalised Poisson distributions truncated to a maximum value of 8 units. We use 10 segments for the piecewise linearisation of the loss functions associated with the stochastic demand.
    \noindent
    The inventory capacity of each retailer is 8 units while the capacity of the vehicle is 10 units. 
    
    The test bed comprises 432 instances, the full factorial of the values considered for the following parameters: road networks, sets of retailer nodes, initial inventory level of retailers, unit penalty costs (in monetary units (MU)) and demand patterns faced by retailers over the planning horizon, as Poisson distribution described by their mean parameters. The possible values considered are described in Table \ref{tab:testbedval}.
    
    For our tests we consider a vehicle with a Gross Vehicle Weight (GVW) of 22 tonnes (t) with a maximum payload of 10t. The capacity of the battery is 150 kWh. The cost of one kWh of energy from ERS or battery is taken as 1 MU ($C^e=1$), while one kWh of fuel costs 3 MU ($C^f=3$). Regarding the battery of the vehicle, we consider 20 possible levels for the state of charge.
    
    As mentioned at the beginning of this Section we must ensure that operations of charge and discharge of battery are closed for the state space considered. To this end we introduce new parameters:
    
    \begin{equation}
        \label{eq:energy_levels}
        \varepsilon^{W_t}_{ij} = \text{round}\left(\left(\alpha_{ij} W_t + \beta_{ij}\right)\frac{20}{150}\right)
    \end{equation}
    
    which account for required energy, expressed in battery level units, of each arc and vehicle weight possible. Equation \ref{eq:unbounded_flow} is changed accordingly:

    \begin{equation}
        \label{eq:unbounded_flow_modified}
        b^u_t = b_{t-1} + \sum_{i=1}^N \sum_{j=1}^N T^{ij}_{t-1} \left(s_{ij} - \varepsilon^{W_{t-1}}_{ij}\right)\\
    \end{equation}


    \subsection{Results and discussion}
    \label{sec:results}
    
    The results of our computational study are summarised in Table \ref{table:results}, which shows the mean, median and standard deviation of the percentage error of the solutions obtained by the MILP heuristic when compared with the optimal solution obtained by the SDP formulation, when pivoting all the parameters considered in our test bed. 
    \noindent
    Overall, the heuristic obtains solutions which costs are 3.72\% above the optimal.
    \noindent
    It is worth noticing that whereas the optimal solution is achieved by reviewing the state of the system on every period, after the demand at the retailers is realised, being able to adapt decisions, our heuristic focuses on an static uncertainty strategy \cite{Bookbinder1988}, providing a fixed policy set at the beginning of the planning horizon. 
    \noindent
    Our future work includes investigating the performance of our heuristic on a receding horizon setting: in this setting the MILP heuristic is used over periods $t,\ldots, T$ (starting with $t=1$) but only decisions on period $t$ are considered. After the demand is observed the state of the system is updated (inventory level of retailers and vehicle load, position of the vehicle on the graph and level of battery of the battery). The same process is successively repeated over periods $t+1,\ldots,T$ until reaching the end of the time horizon. This approach has proven to considerably improve cost performances of lot sizing problem policies such as the one considered in our heuristic \cite{DuralSelcuk2020}.
    
    \section{Case study}
    \label{sec:kent}
    
    \begin{table}[]
    \caption{Summary of the solutions, indicating quantities delivered to retailers and the order of visits to the main nodes: Dover (D), Folkestone (F), Ashford (A), Canterbury (C), Maidstone (M), Sittingbourne (M)}
    \label{table:kent}
    \scalebox{0.73}{
    \begin{tabular}{lllllllllll}
                                                                                   &                                                                                  &                         &                           &                                                                            & \multicolumn{5}{c}{Deliveries to retailers}                                                                           &                                 \\
    \multicolumn{1}{c}{\begin{tabular}[c]{@{}c@{}}Instance \\ number\end{tabular}} & \multicolumn{1}{c}{\begin{tabular}[c]{@{}c@{}}Initial \\ inventory\end{tabular}} & \multicolumn{1}{c}{$p$} & \multicolumn{1}{c}{$C^f$} & \multicolumn{1}{c}{\begin{tabular}[c]{@{}c@{}}Vehicle\\ load\end{tabular}} & \multicolumn{1}{c}{F} & \multicolumn{1}{c}{A} & \multicolumn{1}{c}{C} & \multicolumn{1}{c}{M} & \multicolumn{1}{c}{S} & \multicolumn{1}{c}{Visit order} \\
    1                                                                              & I1                                                                               & 0.1                     & 3                         & 2428.75                                                                    & 854.75                & 828.27                & 745.73                & 0.00                  & 0.00                  & D,F,A,C,D                       \\
    2                                                                              & I1                                                                               & 0.1                     & 6                         & 2417.85                                                                    & 849.16                & 828.27                & 740.42                & 0.00                  & 0.00                  & D,F,A,C,D                       \\
    3                                                                              & I1                                                                               & 0.1                     & 10                        & 0.00                                                                       & 0.00                  & 0.00                  & 0.00                  & 0.00                  & 0.00                  & N/A                             \\
    4                                                                              & I1                                                                               & 0.5                     & 3                         & 3363.47                                                                    & 857.16                & 840.59                & 751.04                & 378.52                & 536.16                & D,F,A,C,S,M                     \\
    5                                                                              & I1                                                                               & 0.5                     & 6                         & 2432.81                                                                    & 642.92                & 823.39                & 844.17                & 122.32                & 0.00                  & D,C,A,F,A,M                     \\
    6                                                                              & I1                                                                               & 0.5                     & 10                        & 2115.73                                                                    & 833.72                & 751.04                & 530.97                & 0.00                  & 0.00                  & D,F,A,C,D                       \\
    7                                                                              & I2                                                                               & 0.1                     & 3                         & 1786.06                                                                    & 556.83                & 537.71                & 691.53                & 0.00                  & 0.00                  & D,F,A,C,D                       \\
    8                                                                              & I2                                                                               & 0.1                     & 6                         & 1802.31                                                                    & 321.57                & 514.75                & 965.99                & 0.00                  & 0.00                  & D,C,C,A,F,D                     \\
    9                                                                              & I2                                                                               & 0.1                     & 10                        & 0.00                                                                       & 0.00                  & 0.00                  & 0.00                  & 0.00                  & 0.00                  & N/A                             \\
    10                                                                             & I2                                                                               & 0.5                     & 3                         & 2969.73                                                                    & 172.92                & 378.52                & 1027.99               & 529.76                & 860.54                & D,C,S,M,A,F                     \\
    11                                                                             & I2                                                                               & 0.5                     & 6                         & 1817.74                                                                    & 430.73                & 526.47                & 860.54                & 0.00                  & 0.00                  & D,C,A,F,D                       \\
    12                                                                             & I2                                                                               & 0.5                     & 10                        & 2160.93                                                                    & 556.83                & 543.12                & 117.39                & 517.98                & 425.61                & D,F,A,M,S,C                    
    \end{tabular}}
    \end{table}
    
\begin{figure}
    \centering
    \caption{Map (a) and graph (b) representation of main HGV roads in Kent}
    \label{fig:mapgraph}
    \begin{subfigure}[b]{0.55\textwidth}
       \includegraphics[width=1\linewidth]{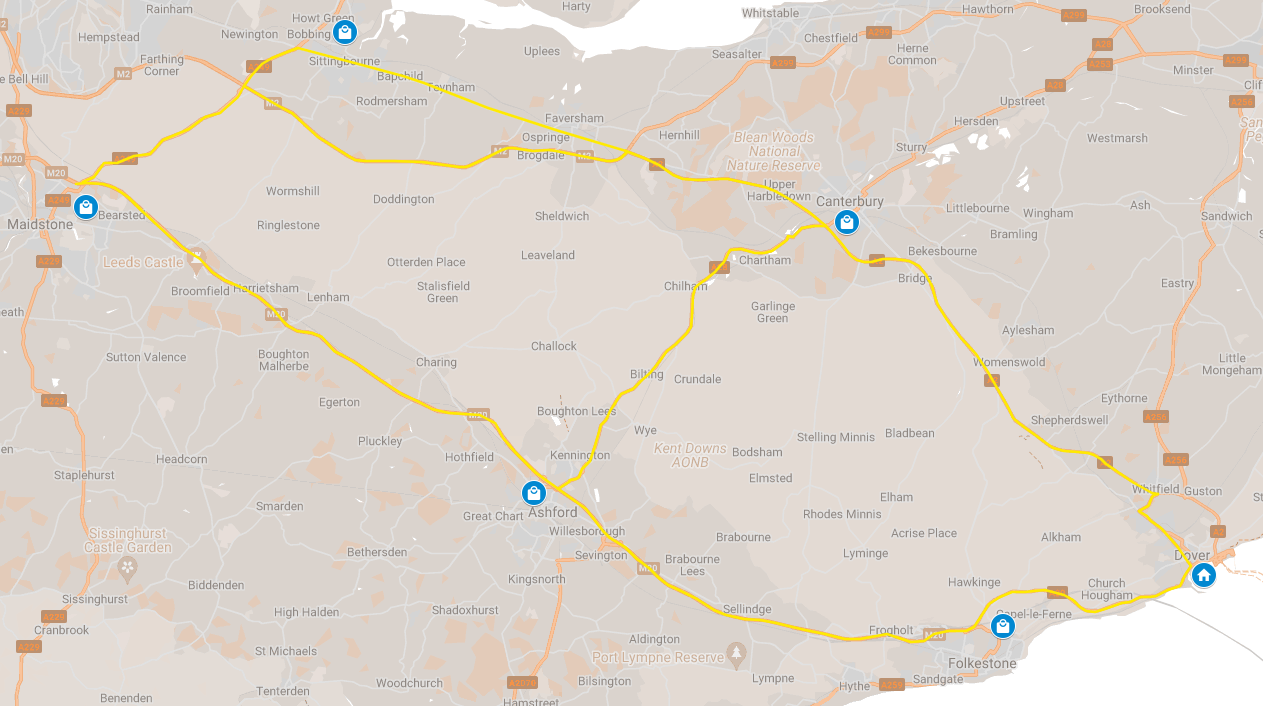}
       \caption{}
       \label{fig:kent_map} 
    \end{subfigure}
    
    \begin{subfigure}[b]{0.65\textwidth}
       \includegraphics[width=1\linewidth]{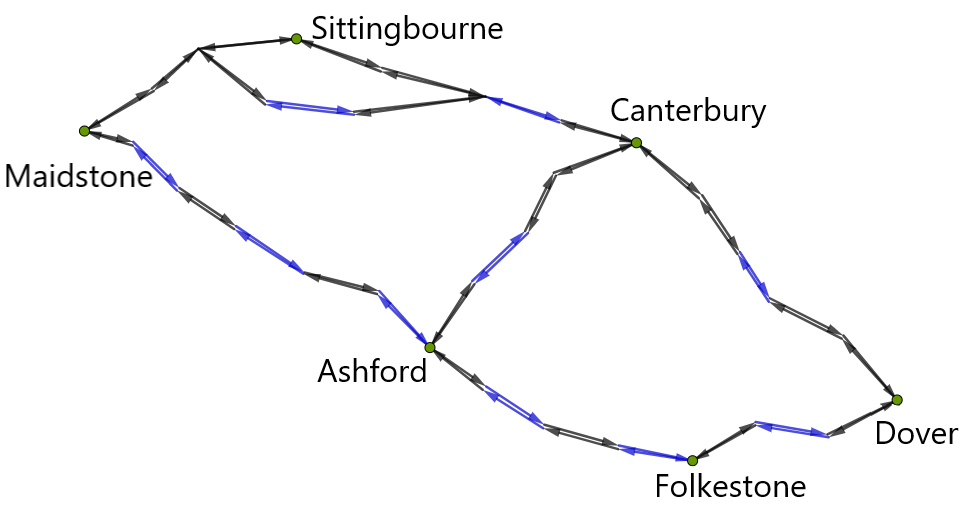}
       \caption{}
       \label{fig:graph}
    \end{subfigure}
    \end{figure}

    In this section we investigate how the relation of the costs of electricity and fuel, and the penalty costs, affect the plan provided by the MILP heuristic. Our new experiments are set up in the region of Kent, in the United Kingdom, creating a realistic instance of the problem, and provide an example of how the delivery quantities and the order of visits to retailers can impact the costs of transportation and reduce the usage of fuel on a ERS network. 
    \noindent
    Figure \ref{fig:mapgraph} shows a map of the region, highlighting the major roads used by HGVs and the graph representation of these roads, in splits of around 5 km between adjacent nodes, making a total of 52 nodes. 
    \noindent
    The depot location is situated at the port of Dover, and we consider five retailer locations situated in densely populated areas of Kent. 
    \noindent
    The electrified stretches of road are marked in blue, and they provide 200 kW of power.
    \noindent
    The HGV considered has a GVW of 17t, with a load capacity of 5t, and it is equipped with a 200 kWh battery.
    \noindent
    Each retailer has a stocking capacity of 2000 kg of product.
    
    The time horizon for our experiments is 25 periods. All five retailers face demand on each period that follow independent normal distributions $\mathcal{N}(50,2)$. We consider two configurations of initial inventory: (I1) all retailers start with 500 kg. of product and (I2) Canterbury and Sittingbourne are out of stock, while the remaining retailers start with 800 kg. The kWh cost of electricity is $C^e=1$, and consider three values for the cost of fuel $C^f = \{3,6,10\}$. We also consider two values for the stockout penalty cost $p = \{0.1,0.5\}$. In all 12 instances the vehicle starts at period 1 at the depot location, with a depleted electric battery and no loaded inventory. We run CPLEX solver with default settings and a time limit of 1 hour. 
    
    Table \ref{table:kent} summarises the solutions reached by the MILP heuristic for the 12 instances, showing the quantities that the vehicle delivers to each retailer and a sketch of the route followed by indicating the order in which the depot and set of retailers are visited. Some trends can be observed in the solutions. First of all, the vehicle load decreases with lower penalty costs and higher fuel costs, to the limiting case of not delivering to any retailers in instances 3 and 9. 
    
    Fuel costs can also have an impact on the routes taken. In I2 instances, Canterbury starts with no inventory and it is usually visited before other retailers (instances 8, 10, 11). However, when fuel costs are high, the vehicle either does not deliver (instance 9), or starts visiting Folkestone and Ashford (instance 12) a route that allows for more battery charging and less fuel usage.

    
    
    \section{Conclusions}
    \label{sec:conclusions}
    In this paper we have dealt with a transportation problem related to emerging technology which allows to charge electric vehicles as they travel, introduced as the Electric Roads Routing Problem, and provided a heuristic for dealing with uncertainty of demand at retailers of the product delivered. The effectiveness of the heuristic was demonstrated by carrying out experiments on an extensive test bed focusing different demand distributions patterns over the time horizon, set of retailer nodes, initial inventory at retailer nodes and on a range of values for the penalty cost. The results show that the heuristic achieves solutions that, overall, are only 3.72\% over the optimal solutions. We tested the heuristic on a realistic case study and showcased how the costs of fuel could impact the routes chosen on a ERS network when delivering to retailers. Future research may investigate the benefit brought by a receding horizon setting and dynamic cut generation strategies as in \cite{Tunc2018} to enhance scalability.
    
	
	\section*{Acknowledgments}
	This project has received funding from the European Union’s Horizon 2020 research and innovation programme under the Marie Skłodowska-Curie grant agreement No 801215, the University of Edinburgh Data-Driven Innovation programme, part of the Edinburgh and South East Scotland City Region Deal and from Costain Group PLC.

	\bibliography{publications}

	\bibliographystyle{plain}
	
	\newpage
	
    \appendix

    	\section{Energy consumption model}
	\label{sec:energymodel}
	
	Travelling costs, either in terms of emissions or economic costs of fuel and electricity usage, are directly related to the energy requirements (kWh) necessary to power a vehicle. To estimate these values we use a simple physics model derived from \cite{Barth2004} and in line with the followed approach by \cite{PRP-2011}.
	
	At any given time, the mechanical power $P$ of a vehicle, as a function of its acceleration $a$ ($m/s^2$) and velocity $v$ ($m/s$) can be estimated as follows:
	\begin{equation}
	    \label{eq:power}
	    P(a,v) = Mav + Mgv \text{sin}\theta + 0.5C_dA\rho v^3 + MgC_r \text{cos}\theta v
	\end{equation}
	
	where $M$ is the mass of the vehicle (kg.), $\theta$ is the road angle in degrees, $C_d$ is a drag coefficient specific for the vehicle, $A$ its frontal area ($m^2$), $g$ is the gravitational constant ($9.81 m/s^2$) and $C_r$ the rolling resistance of the wheels of the vehicle with the road surface.  
	
	For our purposes, we assume the vehicles transit each arc $(i,j)$ with slope $\theta = \theta_{ij}$ cruising at constant speed $v_{ij}$ and therefore, the acceleration is zero, $a=0$. 
	\noindent
	We simply denote the mechanical power on node $(i,j)$ as 
	\begin{equation}
	\label{eq:mechanicalpower}
	    \begin{split}
	        p_{ij} & = P(0,v_{ij})
	    \end{split}
	\end{equation}

	The electric power, and ultimately the battery power that is needed to maintain constant speed $v_{ij}$, depend on multiple factors. As in \cite{Goeke2015} we estimate the battery power required to sustain $v_{ij}$ in arc $(i,j)$ accounting for engine and battery efficiency with a single parameter $\lambda$, obtaining  $P^b_{ij} = \lambda p_{ij}$.
	\noindent
	Finally, the energy required to traverse $(i,j)$ using electric battery, expressed in Joules, can be obtained by multiplying the required power from the battery by the time in seconds spent transiting the arc, $d_{ij}/v_{ij}$, and we express it as a function of the weight of the vehicle $M$:
	
	\begin{equation}
	    \label{eq:battery_energy}
	    \begin{split}
	        r_{ij}(M) = P^b_{ij} (d_{ij}/v_{ij}) & = \lambda (Mg \text{sin}\theta_{ij} + 0.5C_dA\rho v_{ij}^2 + MgC_r \text{cos}\theta_{ij})d_{ij}\\
	        & = \alpha_{ij}M + \beta_{ij}
	    \end{split}
	\end{equation}
	
	where $\alpha_{ij}=\lambda d_{ij}g (\text{sin} \theta_{ij} + C_r \text{cos} \theta_{ij})$ and $\beta_{ij}=\lambda d_{ij}0.5C_d A \rho v^2_{ij}$ are arc constants.
	
	
	In absence of battery energy, the vehicle uses fuel to continue its trip. An estimation of fuel usage and emissions is given in \cite{Barth2004}. Given the engine power output $P$ of the vehicle, the fuel rate is calculated as
	
	\begin{equation}
	    \label{eq:fuelrate}
	    F \approx (kNV + (P/\epsilon + P_a) / \eta)U
	\end{equation}
	
	where $k$ is the engine friction, $N$ is the engine speed, $V$ the engine displacement, $\epsilon$ the drivetrain efficiency of the vehicle, $P_a$ is the amount of engine power related with engine running losses and vehicle accessories such as air conditioner, $\eta$ is a measure of efficiency for diesel engines and $U$ is a factor that depend on $N$ and other constants. 
	
	For the scope of this problem, where travelling speeds are relatively similar on different road segments, we follow a pragmatic approach and approximate fuel directly from the power output $P$, or equally, as a linear transformation of the mechanical energy needed to traverse a specific arc. 

    \section{Symbols used in the MILP formulation}
    \label{sec:MILPsymbols}
    \begin{table*}[ht]
		\begin{align*}
		\intertext{Parameters}
		C                                   &   & \text{Number of retailers}\\
		N                                   &   & \text{Number of nodes}\\
		T                                   &   & \text{Number of time periods}\\
		\delta_{ij}                         &   & \text{Binary parameter indicating if $(i,j)\in \mathcal{A}$}\\
		B                                   &   & \text{Capacity of the battery (kWh) of the vehicle}\\
		p_{ij}                              &   & \text{Mechanical power needed to transit $(i,j)$ at constant speed $v_{ij}$}\\
		s_{ij}                              &   & \text{Supplied energy (kWh) on arc $(i,j)\in \mathcal{A}$}\\
		r_{ij}(M)                              &   & \text{Needed battery (kWh) to traverse $(i,j)\in \mathcal{A}$ with total mass $M$}\\
		C^e                                 &   & \text{Emissions caused per kWh of electric energy used}\\
		C^f                                 &   & \text{Parameter for estimating fuel emissions from required energy}\\
		d^c_t                          &   & \text{Total demand (kg.) on retailer $c \in \mathcal{C}$ in period $t$}\\
		w                            &   & \text{Weight of the vehicle without inventory}\\
		K                            &   & \text{Capacity of the vehicle}\\
		k_c                            &   & \text{Capacity of retailer $c$}\\		
		s_i                           &   & \text{Initial inventory level of retailer $i$}\\
		l                            &   & \text{Initial inventory level of the vehicle}\\
		p                           &   & \text{Penalty costs per unit of lost sales}\\	
		\end{align*}
	\end{table*}

    \begin{table*}[ht]
		\resizebox{\linewidth}{!}{
		\begin{minipage}{\linewidth}
		\begin{align*}
		\intertext{Variables}
		V^i_t					   		    &	& \text{Binary variable set to one iff the vehicle is at node $i$ at time $t$}\\
		T^{ij}_t                            &   & \text{Binary variable set to one iff the vehicle transits from $i$ to $j$ by end of time $t$}\\
		b_t                                 &   & \text{Battery level (in kWh) at time $t$}\\
		b^u_t                               &   & \text{Auxiliary variable: unbounded battery level at time $t$}\\
		b^+_t                               &   & \text{Auxiliary variable: positive battery level at time $t$}\\
		Q^c_t                        &   & \text{Quantity of product (kg.) delivered to retailer $c\in \mathcal{C}$ at time $t$}\\
		S^c_t                        &   & \text{Unsatisfied demand (kg.) on retailer $c\in \mathcal{C}$ at time $t$}\\
		L_t                         &   & \text{Inventory of the vehicle (kg. at time $t$, after replenishment is done}\\
		W_t                          &   & \text{Weight of the vehicle at the end of period $t$}\\
		\end{align*}
		\end{minipage}}
	\end{table*}
\end{document}